\documentclass[10pt]{amsart}
\usepackage{amssymb,amsmath,enumerate,verbatim}%,cases}
\usepackage{eucal,url,stmaryrd
}

\usepackage[pagebackref]{hyperref}
\newtheorem{thrm}{Theorem}[section]
\newtheorem{lemma}[thrm]{Lemma}

\newtheorem{cor}[thrm]{Corollary}

\usepackage[margin=1in]{geometry}
\begin{document}
%\begin{abstract}
%We correct errors in the cited  paper.
%\end{abstract}

% ********************* macroes needed for this paper
%\newcommand{\Om}{\mathbb{R}_{+}^{n-k}\times\mathbb{R}^k}
\newcommand{\Om}{\Omega}
\newcommand{\Rm}{\mathbb{R}^n}
\newcommand{\Rnmk}{\mathbb{R}^{n-k}}
\newcommand{\Rmminusk}{\mathbb{R}^{n-k}}
\newcommand{\Rk}{\mathbb{R}^k}
\newcommand{\Rn}{\mathbb{R}^n}
\newcommand{\Hn}{\mathbb{H}^n}
\newcommand{\bG}{\mathbb{R}^n}

%%%%%%%%%%%%%%%%%%%%%%%%%%%%%%%%%%%%%%%%%%%%%%%%App 2 macroes%%%%%%%%%%%
%\newcommand{\dom}{{\overset{o}{\mathcal{W}}}\,^{1,p}(\Omega)}
\newcommand{\dom}{{\mathcal{D}}\,^{1,p}(\Omega)}
\newcommand{\domG}{{\mathcal{D}}\,^{1,p}(\bG)}
\newcommand{\domO}{{\mathcal{D}}\,^{1,p}(\Omega)}
\newcommand{\domain}[1]{{\mathcal{D}}\,^{1,2}(#1)}
\newcommand{\pn}[1]{\norm{#1}_{\,\domG}}
\newcommand{\pnn}[2]{\norm{#1}_{\,L^{p}(#2)}}
\newcommand{\twost}{\frac {2}{n-2}}
\newcommand{\pst}{\frac {pn}{n-p}}
\newcommand{\psst}{\frac {p(n-s)}{n-p}}
\newcommand{\qn}[1]{\norm{#1}_{\,L^{p\text{*}}(G)}}
\newcommand{\qnn}[2]{\norm{#1}_{\,L^{p\text{*}}(#2)}}
\newcommand{\dl}[1]{\delta_{#1}}
\newcommand{\tr}[1]{\tau_{#1}}
%%%%%%%%%%%%%%%%%%%%%%%%%%%%%%%%%%%%%%%%%%%%%%%%%%

%%%%%%%%%%%%%%%%%from regulchar2-29-00%%%%%%%%%%
\newcommand{\dis}[1]{\displaystyle #1}
%%%%%%%%%%%%%%%%%%%%%%%%%%%%%%%%%%%%%%%%%%%%%%%%%

\newcommand{\algg}{\mathfrak g}
\newcommand{\Lap}{\mathcal{L}}
\newcommand{\lap}{\triangle}
\newcommand{\gv}{\mathfrak v}
\newcommand{\gz}{\mathfrak z}
\newcommand{\norm}[1]{\lVert#1\rVert}
\newcommand{\abs}[1]{\lvert #1 \rvert}
\newcommand{\e}{\textbf {e}}

\newcommand{\st}{\text{*}}
\newcommand{\dxj}[1]{\frac{\partial {#1}}{\partial x_{j}}}
\newcommand{\dt}{\frac{d}{dt}}
\newcommand{\dxa}[1]{\frac {\partial {#1}} {\partial x_{\alpha}} }
\newcommand{\dr}[1]{\frac{\partial {#1}}{\partial r}}
\newcommand{\dyb}[1]{\frac {\partial {#1}} {\partial y_{\beta}} }
\newcommand{\dxap}[1]{\frac {\partial {#1}} {\partial x_{{\alpha'}}} }
\newcommand{\dybp}[1]{\frac {\partial {#1}} {\partial y_{{\beta'}}} }
\newcommand{\dxayb}[1]{\frac {\partial^2 {#1}} {\partial x_{\alpha}{\partial y_{\beta}}} }
\newcommand{\dxaxap}[1]{\frac {\partial^2 {#1}}{\partial x_{\alpha}{\partial y_{{\alpha'}}}}}
\newcommand{\dxaxa}[1]{\frac {\partial^2 {#1}}{\partial x_{\alpha}{\partial y_{\alpha}}}}
\newcommand{\dybybp}[1]{\frac {\partial^2 {#1}}{\partial y_{\beta}{\partial y_{{\beta'}}}}}
\newcommand{\R}{\mathbb{R}}
\newcommand{\Cn}{\mathbb{C}^n}
\newcommand{\dell}{\delta_{\lambda}}
\newcommand{\zeps}{\zeta_{\epsilon}}
\newcommand{\Ab}{\Bar{A}}
\newcommand{\nz}{\lvert z \rvert}
\newcommand{\nw}{\lvert w \rvert}
\newcommand{\zb}{\Bar{z}}
\newcommand{\wb}{\Bar{w}}
\newcommand{\bx}{\boldsymbol x}
\newcommand{\by}{\boldsymbol y}

%%%%%%%%%%%%%%%%%%%%%%%%%%%%%%%%%%%%%%%%%%%%%%%%%
\newcommand{\vx}{\boldsymbol x}
\newcommand{\vy}{\boldsymbol y}
\newcommand{\vz}{\boldsymbol z}
\newcommand{\tphi}{\tilde\phi}
\newcommand{\dx}[1]{\frac {\partial {#1}} {\partial x} }
\newcommand{\dxi}[1]{\frac {\partial {#1}} {\partial \xi} }
\newcommand{\dy}[1]{\frac {\partial {#1}} {\partial y} }
\newcommand{\deta}[1]{\frac {\partial {#1}} {\partial \eta} }

%%%%%%%%%%%%%%%%%%%%%%%%%%%%%%%%%%%%%%%%%%%%%%%%%%
\newcommand{\dmu}[1]{\frac {dz} {\abs x^{#1}} }
\newcommand{\dmus}{\frac {dz} {\abs x^{s}} }

\newcommand{\ueps}{{u_\epsilon}}
\newcommand{\veps}{{v_\epsilon}}
\newcommand{\bueps}{{{\bar u}_\epsilon}}

%%%%%%%%%%%%%%%%%%%%%%%%%%%%%%%%%%%%%%%%%%%%%%%%%%%%%%%%%%%%%%%%%
\title[Corrigenda to "Hardy-Sobolev inequalities.."]{Corrigenda to "$L^p$ estimates and asymptotic behavior for finite energy solutions of extremals
 to Hardy-Sobolev inequalities", Trans. Amer. Math. Soc. 363 (2011), no. 1, 37--62}
 
\date{\today}

\author{Dimiter Vassilev}
\address[Dimiter Vassilev]{ Department of Mathematics and Statistics\\ University of New Mexico\\
Albuquerque, New Mexico, 87131}\email{vassilev@unm.edu}

\maketitle

The  claim and the proof of \cite[Theorem 2.9]{V11} are not correct as stated. We give a  correction in Section \ref{ss:add assumption}. In addition, we correct several typos in the text and strengthen  slightly some results in \cite{V11}.

\section{Corrigenda} Throughout we use the notation of \cite{V11}.

\subsection{}  Except for the  introduction of \cite{V11}, throughout the paper the assumption $0\leq s\leq p$ should be $0\leq s<p$. While the case $s=p$ is included in the Hardy-Sobolev inequality recalled in \cite[Theorem 1.1]{V11}, the problems considered in the rest of the paper assumed tacitly that $s<p$. The latter condition was implied in some, but not all, of the statements in \cite{V11} by some of the made assumptions. 

\subsection{} The following corrections should be made in the statement of \cite[Theorem 2.1)]{V11}.  The assumption on the exponent $t$ should have been $0\leq t<\min \{p,k\}$ rather than the more restrictive $0\leq t<\min \{p,s\}$. The same correction applies to the fourth line on \cite[p.41]{V11}.

\subsection{} In identity \cite[(2.21)]{V11} the norms of $u$ should not be raised to the power $q$.

\subsection{} In the proof of \cite[Theorem 2.4)]{V11}, the formula for $r'$ should be $r'=p^*/(p^*(s)-p)$ not $r'=p^*(s)/(p^*(s)-p)$ , noting that $p^*=p^*(0)$, and we use Sobolev's inequality (i.e., Hardy-Sobolev' inequality with $s=0$) in order to see that $V\in L^{r'}(\Omega)$.

\subsection{}  The exponents $p$ and $p'$ appearing in the denominators of the first line of \cite[(2.51)]{V11} should be switched. Thus, the correct form of \cite[(2.51)]{V11} is
\begin{multline}\label{e:holder}
\abs{\nabla u}^{p-2}\alpha^{p-1} F(u)\nabla u \cdot \nabla \alpha\
\leq\ \frac {\epsilon }{p'}\alpha^p\abs{\nabla u}^p u^{-1} F(u)\ +\
\frac {\epsilon^{-p/p'}}{p}
\abs{\nabla\alpha}^p F(u)u^{p/{p'}}\\
\leq\ \frac {\epsilon }{p'} \alpha^p \abs{\nabla G}^p \ +\
C\frac {\epsilon^{-(p-1)}}{p} q^{p-1} \abs {\nabla \alpha}^p G^p
\end{multline}
with $C=C(p)$. Therefore, the correct form of the inequality following \cite[(2.51)]{V11} is
\[
\text{LHS}\ \geq\ (1-(p-1)\epsilon)\int \alpha^p \abs{\nabla G}^p\, dz \
-\ C\epsilon^{-(p-1)}q^{p-1}\int \abs {\nabla \alpha}^p G^p\, dz
\]
and then we fix $\epsilon=\frac {1}{2(p-1)}$ so that $(1-(p-1)\epsilon)=1/2$.

\subsection{}\label{ss:add assumption} %Additional Assumption in \cite[Theorem 2.9]{V11}

Here, we correct  the statement of \cite[Theorem 2.9]{V11} by adding an additional assumption whose origin can be traced to \cite[Lemma (4.2)]{Bando89} and  then in section \ref{ss:proof of T2.9} we give the necessary modification of the proof. 

\textbf{Additional Assumptions in \cite[Theorem 2.9]{V11}}.
\emph{In addition to the assumptions of \cite[Theorem 2.9]{V11}, suppose that $V\in L^{t_o}$ for some $t_o>r'$ and there are constants $R_o$ and $K_o$  such that we have
\begin{equation}\label{e:decay of V}
\int_{\Rn\setminus B(0,R)}|V|^{t_o}dz \leq \frac {K_o}{R^{(p-s)t_o-n}},
\end{equation}
for  all  $R\geq R_o$.}

It is worth pointing the following example from \cite[Example 3, p. 51]{E}  where it is observed that $$u(z)=|z|^{2-n}\ln |z|$$ has  finite energy outside the unit ball and satisfies the equation
\[
\Delta u=-Vu,\qquad\text{with}\qquad V(z)=(n-2)\frac {1}{|z|^2\ln |z|}.
\]
However, $u$ does not have the fast decay at infinity, i.e., the same as the fundamental solution of the Laplacian. Furthermore, $V\in L^{t_o}$ for any $t_o\geq n/2$.  On the other hand, we note that $V$ does not satisfy \eqref{e:decay of V} since
\[
\int_{|z|>R} V^{t_o} dz\approx \frac {\left( \ln R\right)^{t_o-1}}{(t_o-1)R^{2t_o-n}} +\dots.
\]
Therefore, the claim of \cite[Theorem 2.9]{V11} does not hold with the assumptions of the paper. On the other hand, as pointed above, this counterexample is excluded by the added assumption \eqref{e:decay of V}.

\subsection{Proof of \cite[Theorem 2.9]{V11} with  the extra assumption \eqref{e:decay of V}}\label{ss:proof of T2.9} The argument in the proof of \cite[Theorem 2.9]{V11} is incomplete due to the dependence on $q$ of the constant $C$ in the first integral in the right-hand-side of equation \cite[(2.52)]{V11}, which follows from equation   \cite[(2.16)]{V11}. The correct form of \cite[(2.52)]{V11} is with $C$ replaced by $Cq^{p-1}$ in the first term in the right-hand-side of equation \cite[(2.52)]{V11}, i.e.,
\begin{equation}\label{e:2.52 corrected}
\int_\Omega \abs {\nabla(\alpha G)}^p\, dz\ \leq\ Cq^{p-1}\left (\int_{\Omega} |V|\frac
{\alpha^pG^p}{|x|^s}\,dz \ +\  \int_\Omega \abs {\nabla \alpha}^p
G^p\, dz\right ),
\end{equation}
where $C=C(p)$ is a constant independent of $q$. In particular,
we cannot do infinitely many iterations using \cite[(2.53)]{V11} and obtain the claimed estimate. \footnote{ The author thanks Professor Annunziata Loiudice for pointing the erroneous use of equation \cite[(2.16)]{V11} in the claimed (incorrect) form of \cite[(2.52)]{V11}.}

With the added additional assumptions on $V$ stated in section \ref{ss:add assumption} the proof  follows from the original argument in  \cite[Theorem 2.9]{V11} with the following corrections. The general outline is that to handle  the term involving the potential, which will be kept in the right-hand side, we will rely on the new assumption \eqref{e:decay of V} and several inequalities rather than the original use of \cite[(2.54)]{V11}. Starting with \eqref{e:2.52 corrected} where $\alpha$ is a smooth cut-off function to be specified in a moment,  the  first integral in the right-hand side of \eqref{e:2.52 corrected} can be bounded as follows using H\"older's inequality,
\begin{equation}
\int_{\Omega} |V|\frac
{\alpha^pG^p}{|x|^s}\,dz \leq \left (\int_{\mathrm{supp}\, \alpha} |V|^{t_o}\,dz\right)^{1/{t_o}}\,\left( \int_\Omega \frac {|\alpha G|^{p^*( st)}}{|x|^{st}} \, dz\right)^{1/t}\, \left( \int_\Omega   {|\alpha G|^{p}}\,dz\right)^{1/\kappa},
\end{equation}
where $t$ and $\kappa$ satisfy
\[
\frac {1}{t_o}+\frac {1}{t}+\frac {1}{\kappa}=1 \qquad \text{and}\qquad  \frac { {p^*(st)}}{t}+\frac {p}{\kappa}=p.
\]
The last identities together with the assumptions $t_o>r'=\frac {n}{p-s}$ and  $s<p<n$, imply that we have
\begin{equation}\label{e:exp in holder}
t=\frac {pt_o}{n-p+st_o}>1,\qquad \kappa=\frac{pt_o}{(p-s)t_o-n}>1, \qquad \kappa'=\frac {tp}{p^*(st)}, \qquad \text{and}\qquad 0\leq st<p.
\end{equation}
Therefore, with $\epsilon>0$ the inequality $ab\leq \epsilon \frac
{a^{\kappa'}}{\kappa'}+\epsilon^{-(\kappa-1)}\frac {b^{\kappa}}{\kappa}$, the Hardy-Sobolev inequality and the identity  $\kappa'/t=p/p^*(st)$ give
\begin{multline}\notag
\int_{\Omega} |V|\frac
{\alpha^pG^p}{|x|^s}\,dz \leq \frac {\epsilon}{\kappa'}\left( \int_\Omega \frac {|\alpha G|^{p^*( st)}}{|x|^{st}} \, dz\right)^{\kappa'/t}\ +\  \frac {\epsilon^{-\kappa/\kappa'}}{\kappa}\left (\int_{\mathrm{supp}\, \alpha} |V|^{t_o}\,dz\right)^{\kappa/{t_o}}\,\left( \int_\Omega  {|\alpha G|^{p}}\,dz\right)\\
\leq \frac {\epsilon}{\kappa'}S^{p}_{p,st}\left( \int_\Omega |\nabla(\alpha G)|^{p} \, dz\right)\ +\  \frac {\epsilon^{-(\kappa-1)}}{\kappa}\left (\int_{\mathrm{supp}\, \alpha} |V|^{t_o}\,dz\right)^{\kappa/{t_o}}\,\left( \int_\Omega  {|\alpha G|^{p}}\,dz\right).
\end{multline}
The above inequality and  \eqref{e:2.52 corrected} show that for some constant $C_1=C_1(p,n)$ we have
\begin{equation}\label{e:almost end}
\int_\Omega \abs {\nabla(\alpha G)}^p\, dz\ \leq  C_1{q^{(\kappa-1)(p-1)}}\left (\int_{supp \alpha} |V|^{t_o}\,dz\right)^{\kappa/{t_o}}\,\left( \int_\Omega  {|\alpha G|^{p}}\,dz\right)\ + \ C_1q^{p-1}\left ( \int_\Omega \abs {\nabla \alpha}^p
G^p\, dz\right ),
\end{equation}
with $\epsilon$ chosen so that
\[
C_0q^{p-1}\frac {\epsilon}{\kappa'}S^p_{p,st}=1/2.
\]
With the choice of the cut-off function $\alpha$ made after \cite[(2.57)]{V11} we have that $\alpha\in C^\infty_o(B(z,r))$ with $\alpha\equiv 1$ on $B(z,\rho)$ where $0<\rho<r<R=|z|/2$, so that the gradient satisfies $|\nabla \alpha|\leq 2/(r-\rho)$. This will be used to estimate the second integral in the righthand side of \eqref{e:almost end}.
On the other hand, the new extra assumption \eqref{e:decay of V} shows that for $R\geq R_o$ we have
\begin{equation}\label{e:control norm V}
\left (\int_{\mathrm{supp}\, \alpha} |V|^{t_o}\,dz\right)^{\kappa/{t_o}}\leq \frac {K_o^{\kappa/{t_o}}}{R^p}\leq \frac {K_o^{\kappa/{t_o}}}{(r-\rho)^p},
\end{equation}
recalling the value of $\kappa$ given in \eqref{e:exp in holder}.
From \eqref{e:almost end} in which we use inequality \eqref{e:control norm V}, the gradient estimate for $\alpha$ and $\kappa>1$, we obtain
\begin{equation}\label{e:new2.56}
\int_\Omega \abs {\nabla(\alpha G)}^p\, dz\ \leq C_2 \frac {q^{\kappa(p-1)/q}} {(r-\rho)^{p/q}} \int_{B_r(z)}  {G}^p\, dz,
\end{equation}
with a constant $C_2=C_2(s,p,n, K_o,R_o)$.
Finally, an application of the Hardy-Sobolev inequality and Fatou's theorem to inequality  \eqref{e:new2.56}  gives with $\delta=p^*/p>1$ the inequality
\begin{equation}\label{e:base for Moser infty}
\left(\int_{B_\rho(z)} \abs {u}^{q\delta}\, dz\right)^{1/{\delta q}}\ \leq  C_3 \frac {q^{\kappa(p-1)/q}} {(r-\rho)^{p/q}} \left (\int_{B_r(z)}  \abs {u}^q\, dz\right )^{1/q}
\end{equation}
with a constant $C_3=C_3(s,p,n, K_o,R_o)$. The last inequality is the analogue of inequality \cite[(2.58)]{V11} where, now, we have full control of the constant $C$.
At this point the Moser iteration argument described after \cite[(2.58)]{V11} gives the claimed inequalities \cite[(2.47)]{V11} and \cite[(2.48)]{V11}
noting that since $q_j=q_0\delta^j, j=0,1,\dots,$ we have
\[
p\sum_{j=0}^\infty  {1}/{q_j}=p\frac {p^*}{p^*-p}=\frac {n}{q_o}.%=p\frac {np}{np-p(n-p)}
\]
This completes the (correction) to the proof of \cite[Theorem 2.9]{V11}.

\subsection{Addition to \cite[Section 2]{V11}}
The following Corollary strengthens \cite[Theorem 2.9]{V11} in the case of an exterior domain.
 \begin{cor}\label{c:improved Lp}
  Let all of the assumptions of \cite[Theorem 2.9]{V11} (including the  Additional Assumption \eqref{e:decay of V}) hold. If  the complement of $\Omega$ is a compact subset of $\Rn$, then  the claim of \cite[Theorem 2.9]{V11} is valid for any $q_o >p-1$ rather than $q_o\geq p$.
 \end{cor}

 \begin{proof}
 The proof follows the steps of the above proof of \cite[Theorem 2.9]{V11}, so we only indicate the differences.  Since  $V\in L^{t_o}$ for some $t_o>r'$ from \cite[Theorem 2.1]{V11} we have that $u$ is a bounded function. In particular,  in the proof of \cite[Theorem 2.9]{V11} we can use
   $$G(t)=t^{q/p} \qquad\text{and}\qquad  F(t)=\left(\frac {q}{p} \right)^p\frac {t^{q-p+1}}{q-p+1}$$
   by fixing an $l$ sufficiently large in \cite[(2.15)]{V11} so that $u\leq l$.
   Notice that, now, $F$ and $G$ satisfy $|G(u)|^p=u^q$ and  $F(u)=u|G'(u)|^p$ hence we have again \cite[(2.16)]{V11} and \cite[(2.49)]{V11}. Furthermore, in the case when the complement of $\Omega$ is a compact subset of $\Rn$, the cut-off function $\alpha$  used in the proof of \cite[Theorem 2.9]{V11}, recalled  here above inequality \eqref{e:control norm V},   will be a function with compact support in $\Omega$ for all $R$ sufficiently large.
   Thus, for any $q>p-1$ and $\varepsilon>0$ a positive constant %and $F$ and $G$ defined above
   we have that $\alpha^pF(u+\varepsilon)\in D^{1,p}(\Omega)$ can be taken as a test function in the beginning of the proof of \cite[Theorem 2.9]{V11}. The  proof continues as in section \ref{ss:proof of T2.9} leading  to \eqref{e:new2.56}  in which now we have $G=G(u+\varepsilon)$. Letting $\varepsilon\rightarrow 0$ we obtain \eqref{e:base for Moser infty} for any $q>p-1$ such that $u\in L^q$.
\end{proof}

\subsection{Correction to \cite[Theorem 2.10]{V11}.} First, the statement of \cite[Theorem 2.10]{V11} is about equation \cite[(2.5)]{V11} as stated before the theorem. Thus, in the statement of \cite[Theorem 2.10]{V11} \emph{we assume \cite[(2.24)]{V11} with $V_o\equiv 0$}, i.e.,  $u$ is a weak non-negative solution of
\begin{equation}\label{e:eq in theorem 2.10}
 -\ \text{div}\, (|{\nabla}u|^{p-2}{\nabla}u)\ \leq \  R(z)\frac {\
u^{p-1}}{|x|^s} \qquad \text{in} \quad \Omega
 \end{equation}
with $R\in L^{r'}\cap L^{t_o}$ for some $t_o>r'$. Furthermore, since \cite[Theorem 2.10]{V11} is stated as a direct corollary of  \cite[Theorem 2.1, Theorem 2.5 and Theorem 2.9]{V11} and the latter required the additional assumption as described earlier, we need to \emph{add the extra assumption that the potential $R$ satisfies the condition \eqref{e:decay of V} (with $V$ replaced by $R$)}. Finally, \emph{the domain $\Omega$ should be assumed to have a compact complement} for reasons we explain below.

The condition that $\Omega$ is an exterior domain (the complement is compact) is imposed due to the fact that the decay rate  $u(z)\leq C |z|^{n/q_0}$ given by \cite[Theorem 2.9]{V11} is valid for $q_0\geq p$ while the decay of the fundamental solution corresponds to  the limiting exponent $q_0=p^*/p'=n(p-1)/(n-p)$ which is less than $p$. However, in the case of a domain with a compact complement Corollary \ref{c:improved Lp} can be applied since $p^*/p'>p-1$, hence  $q_0$ can be taken as close to $p^*/p'$ as we wish using \cite[Theorem 2.5]{V11} as in the original argument.

\subsection{Correction and Addition to \cite[Section 3]{V11}}
The statement of \cite[Theorem 3.1]{V11} should be corrected by requiring that\emph{ $\Omega=\mathbb{R}^n$, i.e., we are considering entire solutions}\footnote{The author thanks Professor Annunziata Loiudice for pointing  the fact that the given argument requires an extra assumption on the domain.}.

The correction is needed because the proof of \cite[Theorem 3.1]{V11} uses \cite[(3.2)]{V11}, where $p=2$, which implies \cite[(3.6)]{V11}. However, \cite[(3.2)]{V11} relies on the decay  of the solution implied by \cite[Theorem 2.10]{V11}, thus we need to verify the additional assumption on the potential from section \ref{ss:add assumption}, which in this case is $R\, u^{2^*(s)-2}$ with $R\in L^\infty$. The verification  detailed below uses the scaling invariance of the equation when the domain is the entire space. In fact, we give the following  general result valid not only for solutions to the equations modeled on the scalar curvature equations, but also for \emph{entire} non-negative finite energy solutions of equations modeled on the Euler-Lagrange equation of the extremals of the considered Hardy-Sobolev inequality.

\begin{lemma}\label{l:asympt HS}
Let  $R$ be a bounded function, which is non-negative when $p\not=2$. If $u\in \domG$ is a  non-negative weak solution of the inequality
\begin{equation}\label{e:p-extremals}
 -\ \text{div}\, (|{\nabla}u|^{p-2}{\nabla}u)\ \leq \  R(z)\frac {\
u^{p^*(s)-1}}{|x|^s} \qquad \text{on} \quad \Rn,
 \end{equation}
then  for any $0<\theta<1$ there exists a constant
$C_\theta>0$, such that,
\begin{equation}\label{e:asymp_decay_Laplacian}
u(z)\ \leq \ \frac {C_\theta}{1+\abs {z}^{\theta \frac{n-p}{p-1}}}\,
\norm{u}_{D^{1,p}(\Rn)}.
\end{equation}
\end{lemma}

\begin{proof}
 First, exploiting the scale invariance we will show that $u$ has the slow decay at infinity, see \eqref{e:scale decay} below. Second, we will show that $V(z)=R(z)u^{p^*(s)-p}$ satisfies the additional assumptions listed in section \ref{ss:add assumption}. These two observations suffice for the proof of \eqref{e:asymp_decay_Laplacian}. Indeed, the bound \eqref{e:asymp_decay_Laplacian} follows directly from \cite[Theorem 2.5]{V11}  and Corollary \ref{c:improved Lp}. We only need to notice that we can take $q_o$ in \cite[Theorem 2.9]{V11} as close to $p^*/p'$ as we want, since for $p>1$ we have  $p^*/p'=(np-n)/(n-p)>p-1$. We turn to the detailed proofs of the two claims.

 As mentioned in the preceding paragraph, in the first step of the proof we will show that for all $|z|/2>R_o$ any solution of \eqref{e:p-extremals} has the following "slow decay" property
\begin{equation}\label{e:scale decay}
u(z)\leq \frac {C}{|z|^{(n-p)/p}},
\end{equation}
with a constant $C$ depending on $p$, $n$ and $\norm {u}_{\domO}$, see \cite{Z} for the case of the Yamabe equation.

 In fact, if $u$ is a solution of \eqref{e:p-extremals} then for any $\lambda>0$ the function
\begin{equation}\label{e:scaled function}
v(w)=\lambda^{(n-p)/p}u(\lambda w), \quad z=\lambda w,
\end{equation}
is also a finite energy solution of \eqref{e:p-extremals} with $R(z)$ replaced with $R_\lambda(w)=R(\lambda w)$  \footnote{ The author thanks Professor Annunziata Loiudice for pointing  the fact that the domain invariance requires the new extra assumption $\Omega=\Rn$.}. Let us notice that $$\norm{R_{\lambda}}_{L^\infty(\R^n)}=\norm{R}_{L^\infty(\R^n)}$$ and also the $D^{1,p}\,(\R^n)$ and $L^{p^*}(\R^n)$ norms are invariant under this scaling, $\norm {u}_{\domO}=\norm {v}_{\domO}$.

In order to show the slow decay, it is then enough to show that there exist constants $R_o$ and    $C$, depending only on $p$, $n$, and the invariant under the scaling norms, such that for all $z_o$ with  $\lambda=|z_o|/2>R_o$   we have on the ball $B(w_o,1)$, $w_o=z_o/\lambda$, the estimate
\begin{equation}\label{e:bound for v}
\max_{w\in B(w_o,1)}v(w) \leq C.
\end{equation}
Indeed, \eqref{e:bound for v} implies
\[
\left(\frac {|z_o|}{2}\right)^{(n-p)/p}u(z_o)\leq \left(\frac {|z_o|}{2}\right)^{(n-p)/p}\sup_{|z-z_o|<|z_o|/2 } u( z)=\sup_{\left\vert \lambda w-z_0 \right\vert<\lambda}\lambda^{(n-p)/p} u(\lambda w) =\sup_{\left\vert w-w_0 \right\vert<1}v(w)\leq C,
\]
which gives \eqref{e:scale decay}.

Let us observe that  the desired bound \eqref{e:bound for v} is suggested by the local version of \cite[Theorem 2.1]{V11} stated in \cite[Remark 2.3]{V11}. However,  we have to make sure that the local supremum bound is independent of $\lambda$, therefore we provide the details of the argument.

  Let us show  that  for any fixed $t_o>r'$  we have $V(w)\equiv V_\lambda(w)=R(\lambda w)v^{p^*(s)-p}(w)\in L^{t_o}$ near  the point $w_o=z_o/\lambda$, $\lambda=|z_o|/2$, with a uniform bound on the norm as stated in \eqref{e:V uniform bound} below. For this we follow essentially the argument of the proof  of \cite[Theorem 2.9]{V11} after \cite[(2.53)]{V11} (corrected here by \eqref{e:2.52 corrected}), but now using a smooth bump function $\alpha$  depending only on the distance to the point $w_o$ with support in $B(w_o,\rho_2)$ with $\alpha\equiv 1$ on $B{(w_o,\rho_1)}$ for $1/2<\rho_1<\rho_2<1$; furthermore, we take the  function $G$ as in the proof of \cite[Thorem 2.9]{V11} but defined using $v$ instead of $u$. %

Therefore, from \eqref{e:2.52 corrected}, H\"older' and  Hardy-Sobolev' inequalities ( $pr=p^*(rs)$! ), we
have
\begin{equation}
\norm {\nabla(\alpha G)}^p_{L^p} \ \leq\ C_0 q^{p-1}\norm
{V}_{L^{r'}(\text{supp}\, \alpha)}\,\norm {\nabla(\alpha G)}^p_{L^p}
\ +\  C_0 q^{p-1}\int \abs {\nabla \alpha}^p G^p\, dz.
\end{equation}
Notice that  the function $V$ satisfies
\[
\norm {V}^{r'}_{L^{r'}(supp\, \alpha)}\leq \norm{R}_{L^\infty} \int_{B(w_o,1)} v(w)^{p^*}\,d w= \norm{R}_{L^\infty}\int_{B(z_o,|z_o|/2)} u(z)^{p^*}\,dz\leq \norm{R}_{L^\infty}\int_{\R^n\setminus B(0,|z_o|/2)} u(z)^{p^*}\,dz .
\]
Hence, for any fixed $q_1 >p$, there exists $R_o=R_o(q_1)$ sufficiently large so that
\begin{equation}\label{e:small V}
Cq_1^{p-1}\norm {V}_{L^{r'}(supp\, \alpha)}\leq 1/2
\end{equation}
for $\lambda\equiv |z_o|/2\geq R_o$. Trivially, \eqref{e:small V} holds also for all $q$, such that, $p\leq q\leq q_1$. This shows the validity of \cite[(2.56)]{V11} with the above choice of $\alpha$ and $G$, hence we also have \cite[(2.57)]{V11} with $u$ replaced by $v$ and then \cite[(2.58)]{V11} which in our case reads
\begin{equation}\label{e:moser}
\left(\  \frac {1}{\abs{B(w_o,{\rho_1})}}\int_{B(w_o,{\rho_1})}v^{\delta q} \,
dw\ \right)^{\frac {1}{\delta q}}\ \leq\  \frac
{C^{p/q}q^{(p-1)/q}}{({\rho_2}-{\rho_1})^{p/q}}\ \left(\  \frac
{1}{\abs{B(w_o,{\rho_2})}} \int_{B(w_o,{\rho_2})}v^{ q} \, dw\  \right)^{\frac {1}{
q}},
\end{equation}
where $\delta=p^*/p>1$. After \emph{finitely} many Moser's iterations  of the obtained inequality starting with $q=p^*$ we see that   for any fixed $t_o>r'$  we have $v\in L^{t_o(p^*(s)-p)}$ and the function $V$ satisfies
\begin{equation}\label{e:V uniform bound}
\norm {V}_{L^{t_o}(B_{1/2})}\leq C,
\end{equation}
uniformly for all $\lambda=|z_o|/2\geq R_o$ using again that the $L^{p^*}(\R^n)$ norm is invariant under the scaling.

Next, we follow the proof of \cite[Theorem 2.1b]{V11} with $G=G\circ v$, but in this local case, we use for $1/4<\rho_1<\rho_2<1/2$ a smooth bump function $\alpha$ depending on the distance to $w_o$ such that $\alpha$ is supported in $B(w_o,\rho_2)$ with $\alpha\equiv 1$ on $B(w_o,\rho_1)$. We use again  \eqref{e:2.52 corrected}. The left-hand side is estimates from below  by the Hardy-Sobolev inequality, taking into account the properties of the bump function. 
On the other hand, the first terms in the right-hand side of \eqref{e:2.52 corrected} is estimated from above by H\"older's inequality using \eqref{e:V uniform bound}. The second term is estimated with the help of  the inequalities
$$|\nabla \alpha|\leq 4/(\rho_2-\rho_1)\quad\text{ and }\quad |x|\leq |z|\leq |z-w_o|+|w_o|\leq \rho_2+2=5/2$$
on $B(w_o,\rho_2)\setminus B(w_o,\rho_1)$ which contains the support of $\abs {\nabla \alpha}$. Hence, we obtain
\begin{multline*}
\frac {1}{S^{1/p}_{p,sr}} \left( \int_{B(w_o,\rho_1)}\frac {G^{p^*(s r)}}{|x|^{sr}}\,dz \right)^{p/{p^*(sr)}}
\leq \int_{\R^n} \abs {\nabla(\alpha G)}^p\, dz \leq\ C_0q^{p-1}\left (\int_{\R^n} V\frac
{\alpha^pG^p}{|x|^s}\,dz \ +\  \int_{\R^n} \abs {\nabla \alpha}^p
G^p\, dz\right )\\
\leq  C_1  {q^{p-1}}\,\left [\norm{V}_{L^{t_o}(B(w_o,\rho_2))} \left( \int_{B(w_o,\rho_2)}\frac {G^{p t_o'}}{|x|^{st_o'}}\,dz \right)^{1/{t_o'}} +\left(\frac 25\right)^s \frac {1}{(\rho_2-\rho_1)^p} \int_{B(w_o,\rho_2)}\frac {G^p}{|x|^s}\,dz \right]\\
\leq  C_1(1+\norm{V}_{L^{t_o}(B(w_o,\rho_2))}) \frac {q^{p-1}}{(\rho_2-\rho_1)^p}\,\left [ \left( \int_{B(w_o,\rho_2)}\frac {G^{p t_o'}}{|x|^{st_o'}}\,dz \right)^{1/{t_o'}} +\int_{B(w_o,\rho_2)}\frac {G^p}{|x|^s}\,dz \right]\\
\leq  C_2 \frac {q^{p-1}}{(\rho_2-\rho_1)^p}(1+\rho_2^{n(1-1/t_o')})\,\left( \int_{B(w_o,\rho_2)}\frac {G^{p t_o'}}{|x|^{st_o'}}\,dz \right)^{1/{t_o'}}
\leq C_3 \frac {q^{p-1}}{(\rho_2-\rho_1)^p}\,\left( \int_{B(w_o,\rho_2)}\frac {G^{p t_o'}}{|x|^{st_o'}}\,dz \right)^{1/{t_o'}}
%\left [ \left( \int_{B(w_o,\rho_2)}\frac {G^{p t_o'}}{|x|^{st_o'}}\,dz \right)^{1/{t_o'}} %+\int_{B(w_o,\rho_2)}\frac {G^p}{|x|^s}\,dz \right],
\end{multline*}
after using  also  $4<1/(\rho_2-\rho_1)$, the uniform bound \eqref{e:V uniform bound},  and finally H\"older's inequality taking into account  that $1<t_o'<r$.

 Therefore, using $p^*(rs)=pr$, as noted in \cite[(2.3)]{V11}, and Fatou's lemma show the  existence of a constant $C$, independent of $\lambda=|z_o|/2\geq R_o$, such that,
\begin{equation}\label{e:moser iterate}
\left( \int_{B_{\rho_1}}\frac {|v|^{ qr}}{|x|^{sr}}\,dz \right)^{1/r}
\leq  C \frac {q^{p-1}}{(\rho_2-\rho_1)^p} \,
\left( \int_{B_{\rho_2}}\frac {|u|^{q t_o'}}{|x|^{st_o'}}\,dz \right)^{1/{t_o'}}
\end{equation}
with $0<t_0's<rs$.
An iteration of the above inequality starting with $q=p$, exactly as in the proof of \cite[Theorem 2.1b]{V11} gives the claimed bound \eqref{e:bound for v} using that $rs<k$ by the assumption \cite[(2.4)]{V11}.

 At this point  we can verify the second claim made at the beginning of the proof, i.e.,  that \eqref{e:scale decay} imply the added extra assumptions described in Section \ref{ss:add assumption}.  Indeed, since $V=R(z)u^{p^*(s)-p}\in L^{r'}(\mathbb{R}^n)$,  \cite[Theorem 2.4]{V11} implies that $V\in L^{r'}(\mathbb{R}^n)\cap L^\infty(\mathbb{R}^n)$, while from the slow decay \eqref{e:scale decay}   for $|z|>2R_o$ we have
\[
V(z)\leq \frac {C^{p^*(s)-p}}{|z|^{\frac {n-p}{p}(p^*(s)-p)}}=\frac {C^{p^*(s)-p}}{|z|^{p-s}}
\]
which implies \eqref{e:decay of V}.

This completes the proof of Lemma \ref{l:asympt HS}.
\end{proof}


\begin{thebibliography}{}
\bibitem{E}
Egnell, H.,\  \emph{Asymptotic results for finite energy solutions
of semilinear elliptic equations}.  J. Differential Equations  {\bf
98} (1992),  no.~1, 34--56.

\bibitem{Bando89} Bando, Shigetoshi; Kasue, Atsushi; Nakajima, Hiraku
\emph{On a construction of coordinates at infinity on manifolds with fast curvature decay and maximal volume growth.} Invent. Math. 97 (1989), no. 2, 313--349.

\bibitem{V11} Vassilev, D., \emph{ Lp estimates and asymptotic behavior for finite energy solutions of extremals to Hardy-Sobolev inequalities.} Trans. Amer. Math. Soc. 363 (2011), no. 1, 37--62.

\bibitem{Z}
Zhang, Q.,\  \emph{A Liouville type theorem for some critical
semilinear elliptic equations on noncompact manifolds}, Indiana
Univ. Math. J., {\bf 50} (2001), 1915--1936
\end{thebibliography}
\end{document}